\documentclass[12pt]{article}
\usepackage{amssymb}
\usepackage{amsmath}
\usepackage{latexsym}
\usepackage{mathrsfs}

\numberwithin{equation}{section}
\newtheorem{thm}{Theorem}[section]

\newtheorem{rem}{Remark}[section]
\newtheorem{lem}{Lemma}[section]
\newtheorem{corol}{Corollary}[section]
\newtheorem{example}{Example}[section]
\title{A maximum principle in spectral optimization problems for elliptic operators subject to mass density perturbations\footnote{Published in {\it Eurasian Mathematical Journal}, Volume 4, Number 3 (2013), 70-83.}}

\author{Pier Domenico Lamberti\,\,   and Luigi Provenzano}

\date{\ }

\begin{document}

\newcommand{\rea}{\mathbb{R}}

\maketitle

%
%
%

\noindent
{\bf Abstract:}
We consider eigenvalue  problems
for general elliptic operators of arbitrary order subject to homogeneous  boundary conditions on open subsets of the Euclidean N-dimensional space. We prove stability results for the dependence of the eigenvalues upon variation of the mass density and we prove a maximum principle for extremum problems related to mass density perturbations which preserve the total mass.

\vspace{11pt}


\vspace{6pt}

\section{Introduction}
We consider a general class of elliptic partial differential operators
$$
{\mathcal{L}}u=\sum_{0\le |\alpha|, |\beta |\le m }(-1)^{|\alpha |}D^{\alpha }(A_{\alpha\beta  }D^{\beta }u )
$$
subject to homogeneous boundary conditions on an open subset $\Omega $  of ${\mathbb{R}}^N$ with finite measure.
We assume that the coefficients $A_{\alpha \beta}$ are fixed
bounded real-valued functions such that $A_{\alpha \beta }=A_{\beta \alpha}$ and such that G\r{a}rding's inequality is satisfied. For such operators
we consider the eigenvalue problem
\begin{equation}
\label{classic}
{\mathcal{L}}u=\lambda \rho u \, ,
\end{equation}
where $\rho $ is a positive function bounded away from zero and infinity. Problem (\ref{classic}) admits a divergent
sequence of eigenvalues of finite multiplicity
$$
\lambda_1[\rho]\le \dots \le \lambda_n[\rho ]\le \dots\, .
$$
In this paper we prove a few results concerning  the dependence of $\lambda_n[\rho ]$ upon variation of $\rho$.

Keeping in mind important problems involving harmonic and bi-harmonic operators in linear elasticity (see e.g., Courant and Hilbert~\cite{cohi}), we shall think of the weight $\rho $ as the mass density of the
body $\Omega  $ and we shall refer to the quantity $M=\int_{\Omega }\rho dx$ as the total mass of $\Omega$. In the study of composite materials
it is of interest to know whether it is possible to minimize or maximize the eigenvalues $\lambda_n[\rho]$ under the assumption that the total mass $M$
is fixed (see e.g.,   Chanillo et al.~\cite{chan},  Cox and McLaughlin~\cite{cox, cox1, cox2}, Henrot~\cite{henrot}). In this paper we generalize the results proved in \cite{lamass} for the Dirichlet Laplacian. In particular, we prove the following maximum principle where
we refer to non-zero eigenvalues:\\

{\it All simple eigenvalues and the symmetric functions of multiple eigenvalues of (\ref{classic}) have no points of local maximum or minimum with respect to mass density perturbations preserving the total mass. }\\

See Theorem \ref{nomax} for the precise statement. Moreover, we generalize a result of Cox and McLaughlin~\cite{cox1} and we prove that $\lambda_n[\rho ]$ are weakly* continuous functions of $\rho$, see Theorem~\ref{conti}. This, combined with
the above mentioned principle, implies that if $C$ is a weakly* compact set of mass densities then for non-zero eigenvalues we have: \\

{\it All simple eigenvalues and the symmetric functions of multiple eigenvalues of (\ref{classic}) admit points of maximum and  minimum in $C$ with mass constraint $M={\rm const}$ and such points of maximum and minimum belong to $\partial C$. }\\

See Corollary \ref{cormax} for the precise statement.  The reason why we consider the symmetric functions of multiple eigenvalues and not the eigenvalues themselves is related to well-known bifurcation
phenomena which prevent multiple eigenvalues from being differentiable functions of the parameters involved in the equation.
Moreover, the symmetric functions of multiple eigenvalues appear to be natural objects in the study of extremum problems, see e.g., \cite{buoso, buosoplates, lala2004, lalacri, lamass, lambertistek}. In fact, in this paper
we prove that all
simple eigenvalues and the symmetric functions of multiple eigenvalues are real-analytic functions of $\rho $ and we compute the appropriate formulas
for the Frech\'{e}t differentials which we need for our argument, see Theorem~\ref{sym}.

Theorem \ref{nomax} and Corollary \ref{cormax} are proved for so-called intermediate boundary conditions in which case one of the boundary conditions  is  $u=0$ on $\partial \Omega$  (see condition (\ref{inc}) and Example \ref{exam}). This includes the case of Dirichlet boundary conditions
\begin{equation}\label{dirichlet}
u=\frac{\partial u}{\partial \nu}=\dots =\frac{\partial^{m-1} u}{\partial \nu^{m-1}}=0,\ \ {\rm on }\ \partial \Omega .
\end{equation}
On the other hand, Theorems \ref{conti} and \ref{sym} are proved for a larger class of  homogeneous boundary conditions, including Neumann boundary conditions.  See
Remark~\ref{remneu} for a discussion concerning Neumann-type boundary conditions.

Our work is inspired by the well-known results of  Krein~\cite{krein} and Cox and McLaughlin~\cite{cox, cox1, cox2} concerning the description of optimal mass densities for the eigenvalues of the Dirichlet Laplacian under the additional condition  $A \le \rho \le B$, where $A,B$ are fixed positive constants. Complete solution to this problem for $N=1$ was given in \cite{krein} where explicit fomulas for minimizers and maximizers of all eigenvalues were established. In particular, it turns out that optimal mass densities  are bang-bang solutions, i.e., minimizers and maximizers  satisfy the condition
$(\rho -A )(\rho -B )=0$ on  $ \Omega$. The case $N>1$ is discussed  in \cite{cox1, cox2} where, among other results, it is proved that minimizers and maximizers of the first eigenvalue of the Dirichlet Laplacian are bang-bang solutions. Moreover, Friedland~\cite{friedland} proves that the minimizers of  suitable  functionals of the eigenvalues, in particular of any eigenvalue,  are bang-bang as well. In fact, Friedland~\cite{friedland, friedlandsurv} carries out a deep analysis of extremum problems for the eigenvalues of symmetric compact operators in Hilbert space subject to convex sets of constraints,  which in particular allows to prove that optimal mass densities are bang-bang solutions also for higher order operators with Dirichlet boundary conditions on a bounded open interval (cf. \cite[Thm.~3.3]{friedland}). We mention that explicit solutions for the biharmonic operator in spirit to Krein's results can be found in Banks~\cite{banks1, banks2} and Schwarz~\cite{schwarz}, see  Henrot~\cite[\S~11.4.1]{henrot} for related open problems.

Our approach allows to state a maximum principle concerning all eigenvalues of a quite general class of elliptic operators which can be applied to arbitrary sets $C$ of mass densities, not necessarily convex.

\section{Preliminaries and notation}

Let $\Omega $ be an open set in ${\mathbb{R}}^N$ and  $m\in {\mathbb{N}}$. By $W^{m,2}(\Omega )$ we denote the Sobolev space of
functions in $L^2(\Omega )$ with weak derivatives up to order $m$ in $L^2(\Omega )$, endowed with its standard norm defined by
\begin{equation}\label{sob}
\| u\|_{W^{m,2}(\Omega )}=\bigg(\| u\|^2_{L^{2}(\Omega )}  +  \sum_{|\alpha |= m}\| D^{\alpha }u \|_{L^2(\Omega )}^2\bigg)^{\frac{1}{2}},
\end{equation}
for all $u\in W^{m,2}(\Omega )$.   By $W^{m,2}_0(\Omega )$ we denote the closure in $W^{m,2}(\Omega )$ of the space of $C^{\infty}$-functions with compact support in $\Omega $.

In the sequel, we shall always assume that  $V(\Omega)$ is a fixed closed subspace of $W^{m,2}(\Omega )$ containing $W^{m,2}_0(\Omega )$ and such that the embedding $V(\Omega )\subset L^2(\Omega )$ is compact. Moreover, we shall assume that  $A_{\alpha \beta }\in L^{\infty }(\Omega )$ are fixed coefficients  such that $A_{\alpha \beta }=A_{\beta \alpha }$ for all $\alpha , \beta \in {\mathbb{N}}_0^N$ with $|\alpha |,|\beta |\le m$.

By ${\mathcal{R}}$ we denote the subset of $L^{\infty}(\Omega )$ of those functions $\rho\in L^{\infty }(\Omega )$ such that ${\rm ess }\inf _{\Omega }\rho >0$. Let $\rho \in {\mathcal{R}}$ be fixed.  We consider the following eigenvalue problem
\begin{equation}\label{weak}
\int_{\Omega }\sum_{0\le |\alpha |,|\beta |\le m}A_{\alpha \beta }D^{\alpha }uD^{\beta }\varphi dx=\lambda \int_{\Omega } u\varphi \rho dx,\ \ \forall \varphi \in V(\Omega )\, ,
\end{equation}
in the unknowns $u\in V(\Omega )$ (the eigenfunction) and  $\lambda \in {\mathbb{R}}$ (the eigenvalue). Note that problem (\ref{weak}) is the weak-formulation of problem (\ref{classic}) subject to suitable homogeneous boundary conditions.
The choice of the space $V(\Omega )$ is related to the boundary  conditions in the classical formulation of the problem.
For example, if  $V(\Omega )=W^{m,2}_0(\Omega )$ we obtain  Dirichlet boundary conditions as in (\ref{dirichlet}). If $V(\Omega )=W^{m,2}(\Omega )$ we obtain Neumann boundary conditions.  If
$V(\Omega )=W^{m,2}(\Omega )\cap W^{k,2}_0(\Omega )$, for some $k<m$, we obtain intermediate boundary conditions.  See Example \ref{exam} below. See also Nec\v{a}s~\cite[Chp.1]{nec}.

It is convenient to  denote the left-hand side of equation (\ref{weak}) by ${\mathcal{Q}}[u,\varphi ]$.   It is also convenient to denote by $L^2_{\rho }(\Omega)$ the space $L^2(\Omega )$ endowed with the scalar product defined by
$$
<u_1,u_2>_{\rho }=\int_{\Omega }u_1u_2\rho dx, \ \ \forall \ u_1,u_2\in L^2(\Omega ).
$$
 Note that the corresponding norm $\| u\|_{L^2_{\rho }(\Omega )}$ is equivalent to the standard norm.

We assume that the space $V(\Omega )$ and the coefficients $A_{\alpha\beta }$ are such that G\r{a}r\-ding's inequality holds, i.e., we assume that there exist $a, b>0 $ such that
\begin{equation}\label{gar0}
a\| u\|_{W^{m,2}(\Omega )}^2\le  {\mathcal{Q}}[u,u] +b \|u\|_{L^2(\Omega )}^2\, ,
\end{equation}
for all $u\in V(\Omega )$. Actually, in many cases it will be more convenient to normalize the constants $a,b>0$ in such a way that
\begin{equation}\label{gar}
a\| u\|_{W^{m,2}(\Omega )}^2\le  {\mathcal{Q}}[u,u] +b \|u\|_{L^2_{\rho }(\Omega )}^2\, ,
\end{equation}
for all $u\in V(\Omega )$. For classical conditions on the coefficients $A_{\alpha\beta}$
ensuring the validity of (\ref{gar0}) in the case of Dirichlet boundary conditions we refer to Agmon \cite[Thm.~7.6]{ag}. Moreover, we assume that
there exists $c>0$ such that
\begin{equation}
\label{estder}
{\mathcal{Q}}[u,u]\le c \| u\|_{W^{m,2}(\Omega )}^2,
\end{equation}
for all $u\in V(\Omega )$. Note that since the coefficients $A_{\alpha \beta}$ are bounded,  inequality (\ref{estder}) is always satisfied if $\Omega $ is a bounded open set with Lipschitz boundary (actually, it is sufficient that $\Omega $ is a bounded open set with a quasi-resolved boundary, see Burenkov~\cite[Thm.~6,~p. 160]{bur}).

 Under  assumptions (\ref{gar}), (\ref{estder}), it is easy to prove that problem (\ref{weak}) has a divergent sequence of eigenvalues  bounded below by $-b$. To do so, we consider the bounded linear operator $L$ from $V(\Omega )$ to its dual $V(\Omega )'$ which takes any $u\in V(\Omega )$ to the functional $L[u]$ defined by $L[u][\varphi ]={\mathcal{Q}}[u, \varphi]$, for all $\varphi \in V(\Omega )$. Moreover, we consider
the bounded linear operator $I_{\rho }$ from $L^2_{\rho }(\Omega )$ to   $V(\Omega )'$  which takes any $u\in L^2_{\rho }(\Omega )$ to the functional $I_{\rho }[u]$ defined by
$ I_{\rho}[u][\varphi]=<u, \varphi >_{\rho}$, for all $\varphi \in V(\Omega )$. By inequalities (\ref{gar}), (\ref{estder}) and by the boundedness of the coefficients $A_{\alpha \beta}$, it follows that the quadratic form defined by the right-hand side of (\ref{gar}) induces in $V(\Omega )$ a norm equivalent to the standard norm (\ref{sob}). Hence by the Riesz Theorem, it follows that the operator $L+bI_{\rho }$ is a linear homeomorphism from $V(\Omega )$ onto $V(\Omega )'$. Thus, equation (\ref{weak}) is equivalent to the equation
\begin{equation}
(L+bI_{\rho })^{(-1)}\circ I_{\rho }[u]=\mu u
\end{equation}
where
\begin{equation}\label{rec}
\mu = (\lambda +b  )^{-1}.
\end{equation} Thus, it is natural to consider the operator $T_{\rho} $ from $L^2_{\rho}(\Omega )$ to itself defined by
$$
T_{\rho }:=i\circ (L+bI_{\rho })^{(-1)}\circ I_{\rho },
$$
 where $i$ is the embedding of $V(\Omega )$ into $L^2_{\rho }(\Omega)$. In the sequel, we shall omit $i$ and we shall simply write  $T_{\rho }= (L+bI_{\rho })^{(-1)}\circ I_{\rho }$. Note that
\begin{eqnarray}\lefteqn{
<T_{\rho }u_1,u_2>_{\rho }=I_{\rho}[u_2][ (L+bI_{\rho })^{(-1)}\circ I_{\rho }[u_1] ]}\nonumber \\
& & \quad\quad
=(L+bI_{\rho })[(L+bI_{\rho })^{(-1)}\circ I_{\rho }[u_1]][(L+bI_{\rho })^{(-1)}\circ I_{\rho }[u_2]],
\end{eqnarray}
for all $u_1,u_2\in L^2_{\rho }(\Omega )$. Thus, since the operator $L+bI_{\rho }$ is symmetric it follows that $T_{\rho }$ is a self-adjoint operator in $L^2_{\rho }(\Omega )$. Moreover, if the  embedding $V(\Omega )\subset L^2(\Omega )$ is compact then the  operator $T_{\rho }$ is compact. By inequality (\ref{gar}), $T_{\rho }$ is injective. It follows that the spectrum of $T_{\rho }$ is discrete and consists of a sequence of positive eigenvalues of finite multiplicity converging to zero. Then by  (\ref{rec}) and standard spectral theory, we easily deduce the validity of the following

\begin{lem}\label{gen}
Let $\rho\in {\mathcal{R}}$. Assume that inequalities (\ref{gar}) and (\ref{estder}) are  satisfied for some $a,b,c>0$. Then the eigenvalues of equation (\ref{weak}) have finite multiplicity
and can be represented by means of  a divergent sequence $\lambda_n[\rho ]$, $n\in {\mathbb{N}}$ as follows
\begin{equation}\label{minmax}
\lambda_n[\rho]=\min_{\substack{E\subset V(\Omega )\\ {\rm dim }\, E=n}} \max_{\substack{u\in E\\ u\ne 0}} \frac{\int_{\Omega }\sum_{|\alpha |,|\beta |\le m}A_{\alpha \beta }D^{\alpha }uD^{\beta }u dx}{ \int_{\Omega } u^2\rho dx} \, .
\end{equation}
Each eigenvalue is repeated according to its multiplicity and
\begin{equation}\label{lb}
\lambda_n[\rho ]>-b + \frac{a}{\| \rho \|_{L^{\infty }(\Omega )}}\, ,
\end{equation} for all $n\in {\mathbb{N}}$. Moreover, the sequence $\mu_n[\rho ]=(b+\lambda_n[\rho ])^{-1}$, $n\in {\mathbb{N}}$, represents all eigenvalues of the compact self-adjoint operator $T_{\rho }$.
\end{lem}

\begin{example}  \label{exam}
We consider the case of poly-harmonic operators.
Let $m\in {\mathbb{N}}$. Let $A_{\alpha\beta}= \delta_{\alpha \beta }m!/ \alpha ! $ for all $\alpha, \beta \in {\mathbb{N}}^N$ with $|\alpha |=|\beta |=m$, where $\delta _{\alpha \beta }=1$ if $\alpha =\beta $ and $\delta_{\alpha\beta}=0$ otherwise.
Let $k\in {\mathbb{N}}_0$, $0\le k \le m $ and $V(\Omega )=W^{m,2}(\Omega )\cap W^{k,2}_0(\Omega )$.   Note that (\ref{gar}) and (\ref{estder}) are satisfied for any $b>0$ where $a,c>0$ are suitable constants possibly depending on $b$. Moreover, if $k=m$ and the open set $\Omega $ has finite Lebesgue measure then the embedding $V(\Omega )\subset L^2(\Omega )$ is compact. If $0\le k <m$ and the open set $\Omega $ is bounded and has a Lipschitz continuous boundary then the embedding $V(\Omega )\subset L^2(\Omega )$
is compact (actually it is enough to assume that $\Omega $ is a bounded open set with a quasi-continuous boundary, see Burenkov~\cite[Thm.~8,~p.169]{bur}).  Under these assumptions all corresponding eigenvalues $\lambda_n[\rho]$ are well-defined and non-negative.

Note that if $k=m$ then $V(\Omega )= W^{m,2}_0(\Omega )$ and
 by integrating by parts one can easily realize that the the bilinear form ${\mathcal{Q}}[u,\varphi ]$ can be written in the more familiar form
$$
{\mathcal{Q}}[u,\varphi ] =\left\{\begin{array}{ll}\int_{\Omega }\Delta^{\frac{m}{2}}u\Delta^{\frac{m}{2}}\varphi dx,\ \ &{\rm if }\ m\ {\rm is\ even\, ,}\vspace{2mm}\\
\int_{\Omega }\nabla \Delta^{\frac{m-1}{2}}u\nabla \Delta^{\frac{m-1}{2}}\varphi dx,\ \ &{\rm if }\ m\ {\rm is\ odd\, , }
\end{array}\right.
$$
for all $u,\varphi \in W^{m,2}_0(\Omega )$. In this case we obtain the classic poly-harmonic operator ${\mathcal{L}}=(-\Delta)^m$ subject to
the Dirichlet boundary conditions (\ref{dirichlet}).  Recall that the Dirichlet problem arises in the study of vibrating strings for $N=1$ and $m=1$, membranes for $N=2$ and $m=1$, and  clamped plates for $N=2$ and $m=2$.

In the general case $k\le m$, the classic formulation of the eigenvalue problem   is
$$\left\{\begin{array}{ll}
(-\Delta )^mu=\lambda \rho u,\ \ & {\rm in }\ \Omega,\\
\frac{\partial^j u}{\partial \nu^j}=0, \ \forall\ j=0,\dots , k-1,\ \ & {\rm on }\ \partial \Omega , \\
{\mathcal{B}}_ju=0, \ \forall\ j=1,\dots , m-k, \ \ & {\rm on }\ \partial \Omega ,
\end{array}\right. $$
where ${\mathcal{B}}_j$ are uniquely defined `complementing' boundary operators. See Nec\v{a}s~\cite{nec} for details.    For $N\geq 2$, $m=2$ and $k=1$ we obtain the problem
$$\left\{\begin{array}{ll}
\Delta ^2u=\lambda \rho u,\ \ & {\rm in }\ \Omega,\\
u=0,\ \ & {\rm on }\ \partial \Omega , \\
\Delta u -(N-1)K\frac{\partial u}{\partial \nu}=0,\ \ & {\rm on }\ \partial \Omega ,
\end{array}\right. $$
which is related to the study  of a simply supported plate. Here $K$ is the mean curvature of the boundary of $\Omega$.
See Gazzola, Grunau and Sweers~\cite{gaz} for further details.

Finally, we note that if $m=2$ and $k=0$ then $V(\Omega )=W^{2,2}(\Omega )$ and problem (\ref{weak}) is the weak formulation of a Neumann-type problem
for the biharmonic operator
\begin{equation}
\left\{\begin{array}{ll}
\Delta ^2u=\lambda \rho u,\ \ & {\rm in }\ \Omega,\vspace{1mm}\\
\frac{\partial^2u}{\partial^2\nu}=0,\ \ & {\rm on }\ \partial \Omega ,\vspace{1mm} \\
{\rm div}_{\partial\Omega }[P_{\partial \Omega }[ (D^2u) \nu  ]+ \frac{\partial\Delta u}{\partial \nu }=0,\ \ & {\rm on }\ \partial \Omega ,
\end{array}\right.
\end{equation}
which arises in the study of a vibrating free plate. Here ${\rm div}_{\partial\Omega }$ is the tangential divergence and $P_{\partial \Omega }$ the orthogonal projector onto the tangent hyperplane to $\partial \Omega$. See also Chasman~\cite{chas}.

\end{example}

\section{Continuity and analyticity}

By the min-max principle (\ref{minmax}) it follows that $\lambda_n[\rho ]$ is a locally Lipschitz continuous functions of $\rho\in {\mathcal{R}}$. In fact, one can easily prove that
$$
|\lambda_n[\rho_1]-\lambda_n[\rho_2]|\le \frac{\min\{\lambda_n[\rho_1], \lambda_n[\rho_2]\}+2b}{\min\{{\rm ess}\inf \rho_1 ,{\rm ess}\inf \rho_2\}}\| \rho_1-\rho_2\|_{L^{\infty }(\Omega )}\, ,
$$
for all $\rho_1,\rho_2\in {\mathcal{R}}$  satisfying $\| \rho_1-\rho_2\|_{L^{\infty }(\Omega )} < \min \{{\rm ess}\inf \rho_1 ,{\rm ess}\inf \rho_2\}$. In fact  $\lambda_n[\rho]$ depends with continuity on $\rho$ not only with respect to the strong topology of $L^{\infty }(\Omega )$ but also with respect to the weak* topology, which is clearly more relevant in optimization problems. The following theorem was proved by Cox and McLaughlin~\cite{cox1} in the case of the Dirichlet Laplacian and mass densities uniformly bounded away from zero and infinity. The proof can be easily adapted to the general case. Moreover, it is possible to replace the uniform lower bound for $\rho$ by a weaker assumption.

\begin{thm}\label{conti}
Let $C\subset {\mathcal{R}}$ be a bounded set. Assume that there exist $a,b, c>0$ such that inequalities (\ref{gar}) and (\ref{estder}) are satisfied for all $\rho\in C$.   Then the functions from $C$ to ${\mathbb{R}}$ which take any  $\rho\in C$ to
$\lambda_n[\rho ]$ are weakly* continuous for all $n\in {\mathbb{N}}$.
\end{thm}

{\bf Proof. } Since $C$ is bounded in $L^{\infty }(\Omega )$,
it  suffices to prove that given $\rho \in C$ and  a sequence $\rho_j\in C$, $j\in {\mathbb{N}}$ such that $\rho_j \rightharpoonup^* \rho$ as $j\to \infty $ then   $\lambda_n[\rho_j]\to \lambda_n[\rho]$. To do so, we first prove\footnote{This is clearly trivial if we assume that $0<\alpha \le \rho $ for all $\rho \in C$, in which case $\lambda_n[\rho ]\le \lambda_n[\alpha ]$.} that for each $n\in {\mathbb{N}}$ there exists $L_n>0$ such that $\lambda_n[\rho_j]\le L_n$ for all $j\in {\mathbb{N}}$. Let $n\in {\mathbb{N}}$ be fixed and $u_1, \dots , u_n \in V(\Omega )$ be linearly independent eigenfunctions associated with the eigenvalues $\lambda_1[\rho],\dots , \lambda_n[\rho ]$, normalized by $<u_r,u_s>_{\rho}=\delta_{rs}$ for all
$r,s=1,\dots , n$. Note that
$$
\lim_{j\to \infty }\int_{\Omega }u_ru_s\rho_jdx = \int_{\Omega }u_ru_s\rho dx ,
$$
for all $r,s=1,\dots , n$. Thus
\begin{equation}\label{uni}
\lim_{j\to \infty }\int_{\Omega }\bigg(\sum_{r=1}^n \gamma_ru_r\bigg)^2\rho_jdx = \int_{\Omega }\bigg(\sum_{r=1}^n \gamma_ru_r\bigg)^2\rho dx,
\end{equation}
uniformly with respect to $\gamma =(\gamma_1, \dots , \gamma _n)\in {\mathbb{R}}^n$ with $|\gamma |\le 1$. Let $E$ be the linear space generated by $u_1, \dots , u_n$. By (\ref{uni}) it follows that for any $\epsilon >0$ there exists $j_{\epsilon }\in {\mathbb{N}}$ such that
\begin{eqnarray}\label{unib}\lefteqn{
\frac{\int_{\Omega }\sum_{|\alpha |,|\beta |\le m}A_{\alpha \beta }D^{\alpha }uD^{\beta }u dx}{ \int_{\Omega } u^2\rho_j dx}\le
\frac{\int_{\Omega }\sum_{|\alpha |,|\beta |\le m}A_{\alpha \beta }D^{\alpha }uD^{\beta }u dx}{ \int_{\Omega } u^2\rho dx}}\nonumber \\
& & \qquad\qquad\qquad\qquad\qquad\quad\quad+\epsilon (\lambda_n[\rho ]+2b) \le
\lambda_n[\rho]+\epsilon (\lambda_n[\rho ]+2b)\,
\end{eqnarray}
for all $u\in E$, $j\geq j_{\epsilon}$. By combining (\ref{minmax}) and (\ref{unib}) we deduce that $\lambda_n[\rho_j]\le \lambda_n[\rho ]+\epsilon
(\lambda_n[\rho ]+2b)$
for all $j\geq j_{\epsilon }$, which implies the existence of a uniform bound $L_n$ as claimed  above. The rest of the proof follows the lines of Cox~\cite{cox1}.
Let $u_n[\rho_j]$, $n\in {\mathbb{N}}$ be a sequence of eigenfunctions associated with the eigenvalues $\lambda_n[\rho_j]$ normalized by
$<u_n[\rho_j], u_l[\rho_j]>_{\rho_j}= \delta_{nl}$ for all $n,l\in {\mathbb{N}}$. Note that  ${\mathcal{Q}}[u_n[\rho_j],u_n[\rho_j ]] =\lambda_n[\rho_j]$  for all $j\in {\mathbb{N}}$. By inequality (\ref{gar}),
the sequence $u_n[\rho_j]$, $j\in {\mathbb{N}}$ is bounded in the space $V(\Omega )$ equipped with the norm (\ref{sob}). It follows that possibly passing to subsequences, there exists $\bar u_n\in V(\Omega )$ such that $u_n[\rho_j]$ weakly converges to $\bar u_n$ as $j\to \infty $ in $V(\Omega )$, and there exists   $\bar\lambda_n\in {\mathbb{R}}$ such
$\lambda_n[\rho_j]$ converges to $\bar \lambda_n$ as $j\to \infty $.
Moreover, since the embedding $V(\Omega )\subset L^2(\Omega )$
is compact  we can directly assume that $u_n[\rho_j]$ converges to $\bar u_n$ strongly in $L^2(\Omega )$ as $j\to \infty$. By passing to the limit in the weak equation
$$
{\mathcal{Q}}[u_n[\rho_j], \varphi]=\lambda_n[\rho_j]<u_n[\rho_j],\varphi >_{\rho_j},\ \ \ \forall\ \varphi \in V(\Omega )\, ,
$$
it follows that $\bar \lambda_n$ is an eigenvalue and of problem (\ref{weak}) and $\bar u_n$ a corresponding eigenfunction.
Note that $<\bar u_n,\bar u_l>_{\rho }=\delta_{nl}$ for all $n,l\in {\mathbb{N}}$, hence  $\lambda_n$, $n\in {\mathbb{N}}$ is a divergent sequence.  It remains to prove that
$\bar \lambda_n=\lambda_n[\rho ]$ for all $n\in {\mathbb{N}}$. To do so, assume by contradiction that there exists an eigenfunction $\bar u\in V(\Omega)$ associated with an eigenvalue $\bar \lambda $ of the weak problem (\ref{weak}) such that $<\bar u, \bar u_n >_{\rho }=0$ for all
$n\in {\mathbb{N}}$. Assume that $\bar u$ is normalized by $\| \bar u \|_{\rho }= 1/(b+\bar \lambda )$. By the Auchmuty principle \cite{auch} applied to the operator $L+bI_{\rho}$, we have
\begin{equation}\label{auc}
-\frac{1}{2(b+\lambda_n[\rho_j])}\le  \frac{{\mathcal{Q}}[u,u]+b\|u \|^2_{L^2_{\rho_j }(\Omega )}}{2}-\| u-P_{n-1, \rho_j }u  \|_{L^2_{\rho_j}(\Omega )}\ ,
\end{equation}
for all $u\in V(\Omega )$ and  $n,j\in {\mathbb{N}} $. Here $P_{n-1,\rho_j}u$ denotes the orthogonal projection in $L^2_{\rho_j}(\Omega )$ of $u$ onto the
space generated by $u_1[\rho_j], \dots , u_{n-1}[\rho_j]$ for all $n\geq 2$ and $P_{0,\rho_j}u\equiv 0$. By setting $u=\bar u$ and passing to the limit in (\ref{auc}) as $j\to \infty$, we obtain
$$
-\frac{1}{2(b+\bar \lambda_n)}\le \frac{{\mathcal{Q}}[\bar u,\bar u]+b\|\bar u \|^2_{L^2_{\rho }(\Omega )}}{2}-\| \bar u\|_{L^2_{\rho }(\Omega )} =-\frac{1}{2(b+\bar \lambda )}
$$
for all $j\in {\mathbb{N}}$, which contradicts the fact that $\bar\lambda_n\to \infty $ as $n\to \infty $. \hfill $\Box$\\

By classical results in perturbation theory, one can prove that $\lambda_n[\rho ]$ depends real-analytically on $\rho$ as long as $\rho $ is such that
$\lambda_n[\rho ]$ is a simple eigenvalue. This is no longer true if the multiplicity of $\lambda_n[\rho]$  varies. In the case of multiple
eigenvalues, analyticity can be proved for the symmetric functions of the eigenvalues. Namely, given a finite set of indexes  $F\subset {\mathbb{N}}$, we set
$$
{\mathcal { R}}[F]\equiv \left\{\rho\in {\mathcal { R}}:\
\lambda_j[\rho ]\ne \lambda_l[\rho ],\ \forall\  j\in F,\,   l\in \mathbb{N}\setminus F
\right\}
$$

and

\begin{equation}
\label{sym1}
\Lambda_{F,h}[\rho ]=\sum_{ \substack{ j_1,\dots ,j_h\in F\\ j_1<\dots <j_h} }
\lambda_{j_1}[\rho ]\cdots \lambda_{j_h}[\rho ],\ \ \ h=1,\dots , |F| .
\end{equation}

Moreover, in order to compute formulas for the Frech\'{e}t differentials, it is also convenient to set
$$
\Theta [F]\equiv \left\{\rho\in {\mathcal { R}}[F]:\ \lambda_{j_1}[\rho ]
=\lambda_{j_2}[\rho ],\, \
\forall\ j_1,j_2\in F  \right\} .
$$

Then we have the following result

\begin{thm}
\label{sym}
Assume that there exist $a,b,c>0$ such that inequalities (\ref{gar0}) and (\ref{estder}) are satisfied.
Let $F$ be a finite subset of ${\mathbb{N}}$. Then  ${\mathcal { R}}[F]$ is an open set in $L^{\infty}(\Omega )$ and the functions
$\Lambda_{F,h}$
are real-analytic in ${\mathcal { R}}[F]$. Moreover, if $F=\cup_{k=1}^nF_k$ and
$\rho\in \cap_{k=1}^n\Theta [F_k]$ is such that for each $k=1,\dots , n$ the eigenvalues $\lambda_j[\rho ]$ assume the common value $\lambda_{F_k}[\rho ]$
for all $j\in F_k$, then the differentials of the functions $\Lambda_{F,h}$
at the point $\rho$
are given
by the formula
\begin{equation}
\label{sym2}
d\Lambda_{F,h}[\rho][\dot{\rho}] =-\sum_{k=1}^n
c_k
\sum_{l\in F_k}
\int_{\Omega}u_l^2\dot{\rho}dx\, ,
\end{equation}
for all $\dot{\rho}\in L^{\infty }(\Omega)$, where
$$c_k=
\sum_{\substack{0\le h_1\le |F_1|\\ \dots\dots \\ 0\le h_n\le  |F_n|\\ h_1+\dots +h_n=h }}
{ |F_k|-1 \choose h_k-1 }\lambda_{F_k}^{h_k}[\rho]  \prod_{\substack{j=1\\ j\ne k}}^n { |F_j| \choose h_j }\lambda_{F_j}^{h_j}[\rho] ,
$$
and for each $k=1,\dots , n$, $\{ u_l\}_{l\in F_k}$ is an orthonormal basis in $L^2_{\rho }(\Omega)$ of the eigenspace associated with  $\lambda_{F_k}[\rho ]$.
\end{thm}

{\bf Proof.} We set
$$\tilde \Lambda_{F,h}[\rho ]= \sum_{ \substack{ j_1,\dots ,j_h\in F\\ j_1<\dots <j_h} }
(\lambda_{j_1}[\rho ]+b)\cdots (\lambda_{j_h}[\rho ]+b) \, ,$$
for all $\rho \in {\mathcal{R}}[F]$. Note that by elementary combinatorics, we have

\begin{equation}\label{comb}
\Lambda_{F,h}[\rho]=\sum_{k=0}^h(-b)^{h-k}{
|F|-k
\choose
h-k
}\tilde\Lambda_{F,k}[\rho ]\, ,
\end{equation}
where we have set  $\Lambda_{F,0}=\tilde \Lambda_{F,0}=1$.

By adapting to the operator $L+bI_{\rho}$ the same argument used in \cite{lamass} for the Dirichlet Laplacian, one can prove that
${\mathcal { R}}[F]$ is an open set in $L^{\infty }(\Omega)$ and that $ \tilde \Lambda_{F,h}[\rho ]$ depends real-analytically on $\rho\in {\mathcal { R}}[F]$. Thus, by (\ref{comb}) we deduce the real-analyticity of the functions $\Lambda_{F,h}$.

We now prove formula (\ref{sym2}).  First we assume that $n=1$, hence $F=F_1$ and   $\rho \in \Theta [F]$. For simplicity, we write $\lambda_F[\rho ]$ rather than $\lambda_{F_1}[\rho ]$. The same computations used in \cite{lamass} yields the following formula for the Frech\'{e}t differential $d\tilde\Lambda_{F,h}[\rho]$ of
$\tilde\Lambda_{F,h}$ at the point $\rho\in {\mathcal { R}}[F]$:
\begin{equation}
d\tilde\Lambda_{F,h}[\rho][\dot{\rho}] =
-(\lambda_{F}[\rho ]+b)^{h+1}
{
|F|-1
\choose
h-1
}
\sum_{l\in F}
<dT_{\rho}[\dot \rho][u_l],u_l  >_{\rho }
\, , \ \ \ \ \forall \dot\rho\in L^{\infty}(\Omega ).
\end{equation}

By standard calculus and by recalling that $T_{\rho }u_l=(\lambda_F[\rho]+b)^{-1}u_l$ for all $l\in F$, we have
\begin{eqnarray}\lefteqn{
<dT_{\rho}[\dot \rho][u_l],u_l  >_{\rho }=-b<(L+bI_{\rho })^{-1}dI_{\rho}[\dot \rho](L+bI_{\rho })^{-1}I_{\rho}u_l,u_l  >_{\rho}}\nonumber  \\ & &
+<(L+bI_{\rho })^{-1}dI_{\rho}[\dot \rho]u_l,u_l>_{\rho}=\frac{\lambda_F[\rho]}{\lambda_F[\rho]+b} <(L+bI_{\rho })^{-1}dI_{\rho}[\dot \rho]u_l,u_l>_{\rho}\nonumber  \\
& & =\frac{\lambda_F[\rho]}{(\lambda_F[\rho]+b)^2}\int_{\Omega}u_l^2\dot \rho dx
\end{eqnarray}
hence

\begin{equation}\label{tildif}
d\tilde\Lambda_{F,h}[\rho][\dot{\rho}] =-\lambda_F[\rho ]
(\lambda_{F}[\rho ]+b)^{h-1}
{
|F|-1
\choose
h-1
}
\sum_{l\in F}
\int_{\Omega}u_l^2\dot \rho dx
\, ,
\end{equation}
for all $\dot\rho\in L^{\infty}(\Omega )$. By (\ref{comb}) and (\ref{tildif}) we get
\begin{eqnarray*}\lefteqn{
d \Lambda_{F,h}[\rho][\dot{\rho}]}\\ & &
=-\sum_{k=1}^h\lambda_F[\rho](\lambda_F[\rho]+b)^{k-1}(-b)^{h-k}{
|F|-1
\choose
k-1
}{
|F|-k
\choose
h-k
}\sum_{l\in F}
\int_{\Omega}u_l^2\dot \rho dx \\ & &
=-\lambda_F[\rho]{
|F|-1
\choose
h-1
}\sum_{k=0}^{h-1}{
h-1
\choose
k
}(\lambda_F[\rho]+b)^k(-b)^{h-1-k}\sum_{l\in F}
\int_{\Omega}u_l^2\dot \rho dx,
\end{eqnarray*}
which immediately implies (\ref{sym2}) for $n=1$. We now consider the case $n>1$.  By means of a continuity argument, one can easily see that
there exists
an open neighborhood ${\mathcal{W}}$ of $ \rho$ in ${\mathcal{R}}[F]$ such that ${\mathcal{W}}\subset \cap _{k=1}^n {\mathcal{R}}[F_k]$.
Thus,
\begin{equation}
\label{blocks}
\Lambda _{F,h}= \sum_{\substack{0\le h_1\le  |F_1|, \dots , 0\le  h_n\le |F_n |\\ h_1+\dots +h_n=h}}\,\, \prod_{k=1}^n\Lambda_{F_k,h_k}
\end{equation}
on $ {\mathcal{W}}$. By differentiating equality (\ref{blocks}) at the point $ \rho$ and applying formula (\ref{sym2}) for $n=1$ to each function $\Lambda_{F_k,h_k}$, we deduce the validity of formula (\ref{sym2}) for arbitrary values of $n\in {\mathbb{N}}$. \hfill $\Box$

\section{Maximum principle}

In this section we consider the case of general intermediate boundary conditions. This means that we assume that $V(\Omega )$ is a closed subspace of $W^{m,2}(\Omega )$ satisfying the inclusion
\begin{equation}\label{inc}
V(\Omega )\subset W^{1,2}_0(\Omega )\, .
\end{equation}

Assume that $\Omega $ has finite measure. For all $M>0$ we set
\begin{equation}
L_M=\left\{\rho\in L^{\infty }(\Omega ):\ \int_{\Omega }\rho  dx=M  \right\}
\end{equation}

The following theorem is a generalization of \cite[Thm.~4.4]{lamass} to the case of intermediate boundary conditions.

\begin{thm}\label{nomax}Let  all assumptions of Theorem \ref{sym} hold. Assume in addition  that $\Omega $ has finite measure and inclusion (\ref{inc}) holds.
Then for all $h=1, \dots , |F| $
the map $\Lambda_{F, h} $  of  ${\mathcal{R}}[F] \cap L_M$  to $\mathbb{R}$ which takes
any $\rho\in {\mathcal{R}}[F]\cap L_M$ to   $\Lambda_{F, h}[\rho ] $
has no points of local maximum or minimum $\tilde \rho $ such  that $\lambda_j[\tilde \rho ]$  have the same sign  and  $\lambda_j[\tilde \rho ]\ne 0$ for all $j\in F$.
\end{thm}

{\bf Proof.} It is convenient to consider the real-valued function $M$ defined on $L^{\infty }(\Omega )$ by $M[\rho]=\int_{\Omega }\rho dx$ for all  $\rho \in L^{\infty }(\Omega )$. Assume by contradiction the existence of $\tilde \rho$ as in the statement. Then $\tilde \rho$ is a critical point for the function
$\Lambda_{F,h}$ subject to the mass constraint $M[\rho ]=M$.
This implies the existence of a Lagrange multiplier which means that there exists
$c\in {\mathbb{R}}$ such that  $d\Lambda_{F,h} [\tilde \rho ]=c dM[\tilde \rho ]$ (see  e.g., Deimling~\cite[Thm.~26.1]{dei}).
By formula (\ref{sym2}), it follows that
$$
\int_{\Omega }\left( \sum_{k=1}^nc_k\sum_{l\in F_k}u_l^2 \right) \dot\rho dx
=   c  \int_{\Omega}\dot{\rho}dx,
$$
for all $\dot{\rho}\in L^{\infty }(\Omega )$.  Note that   $c_k$ are non-zero real numbers of the same sign.
Since  $\dot{\rho}$ is arbitrary, it follows that
\begin{equation}\label{lincomb}
\left( \sum_{k=1}^nc_k\sum_{l\in F_k}u_l^2 \right) =c, \ \ {\rm a.e.\ in }\ \Omega .
\end{equation}
Since $u_l\in W^{1,2}_0(\Omega )$, then by a standard argument  one can  prove that the function $(\sum_{k=1}^n\sum_{l\in F_k}(\sqrt{|c_k|}u_l)^2 )^{1/2}$ belongs to the space
$ W^{1,2}_0(\Omega ) $ and
equals $\sqrt{|c|}$ almost everywhere in $\Omega$. As is well-known  the space  $W^{1,2}_0(\Omega )$ does not contain constant functions apart from the function identically equal to zero. Thus $c=0$  and accordingly $u_l=0$ for all $l\in F$, a contradiction.
\hfill $\Box$

\begin{rem}Theorem \ref{nomax} concerns mass densities $\tilde \rho $ such that $\lambda_j[\tilde \rho ]$ do not vanish and  have the same sign for all $j\in F$. This assumption is clearly guaranteed for positively defined operators. Moreover, we note that the sign of the eigenvalues is preserved by small perturbations of $\rho$. Hence our assumption is not much restrictive in the analysis of bifurcation phenomena associated with multiple eigenvalues different from zero.
\end{rem}

Finally, by Theorems \ref{conti} and  \ref{nomax} we deduce the following

\begin{corol}\label{cormax}Let  all assumptions of Theorem \ref{nomax} hold.
Let $C\subset {\mathcal{R}}[F]$ be a $weakly^*$ compact set in $L^{\infty }(\Omega )$. Assume that there exist $a,b>0$ such that
inequality (\ref{gar}) is satisfied for all $\rho \in C$. Let $M>0$ be such that $C\cap L_M$ is not empty.
Assume that the eigenvalues $\lambda_j[\rho ]$ have the same sign and do not vanish for all $j\in F$, $\rho \in C$.
Then for all $h\in \{1,\dots , |F|\}$ the map $\Lambda_{F,h}$ from  $C\cap L_M$ to ${\mathbb{R}}$  which takes
$\rho\in C\cap L_M$ to   $\Lambda_{F, h}[\rho ] $ admits points of  maximum and minimum and all such points belong to $\partial C\cap L_M$.
\end{corol}

{\bf Proof.} Recall that weakly* compact sets are bounded. Thus, by Theorem \ref{conti} the functions $\Lambda_{F,h}$ are weakly* continuous on $C$  hence they admit both maximum and minimum on the weakly* compact subset $ C\cap L_M$ of $C$. By Corollary~\ref{nomax} the corresponding  points of maximum and minimum  cannot be interior points of $C$, hence they belong to $\partial C\cap L_M$. \hfill $\Box$\vspace{12pt}

\begin{example} Consider  the poly-harmonic operators subject to Dirichlet or intermediate boundary conditions as described in Example \ref{exam}.
Let $A,B\in L^{\infty }(\Omega )$ be functions satisfying the condition
$$
0< {\rm ess}\inf_{x\in \Omega } A(x) < {\rm ess}\sup_{x\in \Omega }B(x)<\infty  .
$$
Let $C=\{\rho \in L^{\infty }(\Omega ):\ A\le \rho \le B\}$. Clearly, $C$ is a weakly* compact set. Moreover, since all mass densities $\rho$ are
uniformly bounded away from zero and infinity, inequality
(\ref{gar}) is satisfied for suitable constants $a,b>0$ not depending on $\rho \in C$. Thus Corollary \ref{cormax} is applicable to all non-zero eigenvalues. It turns out that
points of maximum and minimum $\tilde \rho $ should coincide with $A(x)$ or $B(x)$ in a set of positive measure.
\end{example}

\begin{rem}\label{remneu} Condition (\ref{inc}) was used only to guarantee that $V(\Omega )\setminus \{0\}$ does not contain constant functions.
Thus, one may replace condition (\ref{inc}) by slightly more general conditions. For example one may assume that $V(\Omega )\subset W^{1,2}_{0,\Gamma }(\Omega )$ where
$W^{1,2}_{0,\Gamma }(\Omega )$ is the closure in $W^{1,2}(\Omega )$ of $C^{\infty }$-functions vanishing in an open neighborhood of a suitable subset of $\Gamma$ of $\partial \Omega $. In this case, one would talk about mixed-intermediate boundary conditions.

If $V(\Omega )$ is a closed subspace of $W^{m,2}(\Omega )$ containing constant functions different from zero, then
we could argue as in the proof on Theorem \ref{nomax} up to condition (\ref{lincomb}). Thus, in the general case one could simply characterize the critical mass densities of the functions $\Lambda_{F,h}$ as those mass densities for which condition (\ref{lincomb}) is satisfied.
Clearly, in the case of simple eigenvalues condition (\ref{lincomb}) reduces to $u= {\rm const}$ in $\Omega $ which implies that $\lambda =0$. Thus, we conclude that the maximum principle stated in the introduction holds for all simple eigenvalues and all homogeneous boundary conditions under consideration. As for multiple eigenvalues we note that the analysis of condition (\ref{lincomb})  is not straightforward as it may appear
at a first glance. Under suitable regularity assumptions on the eigenfunctions  $u_1$, $u_2$ associated with a double  eigenvalue $\lambda $ of the Neumann Laplacian, one may prove that the condition $u_1^2+u_2^2={\rm const}$ in $\Omega $ implies that $\lambda =0$. However, we do not include such arguments here since we plan to perform a deeper analysis of Neumann and other boundary conditions in a forthcoming paper.
\end{rem}

{\bf Acknowledgments}: We acknowledge financial support  by the research pro\-ject  ``Singular perturbation problems for differential operators", Progetto di Ateneo
of the University of Padova. The first author acknowledges financial support from the research project PRIN 2008 ``Aspetti geometrici delle equazioni alle derivate parziali e questioni connesse''. The second author was partially supported by the Grant MTM2008-03541 of
the MICINN, Spain, and the ERC Advanced Grant NUMERIWAVES.
\\

\noindent {\small
Pier Domenico Lamberti and Luigi Provenzano\\
Dipartimento di Matematica\\
Universit\`{a} degli Studi di Padova\\
Via Trieste,  63\\
35126 Padova\\
Italy\\
e-mail:	lamberti@math.unipd.it \\
e-mail: proz@math.unipd.it

}
\end{document}